\documentclass[11pt,reqno]{amsart}
\usepackage{amsmath,amssymb,amsfonts,mathrsfs}
\usepackage[margin=1.05in]{geometry}
\usepackage{microtype}
\usepackage{booktabs}
\usepackage[colorlinks=true,linkcolor=blue,citecolor=blue,urlcolor=blue]{hyperref}

\newtheorem{theorem}{Theorem}[section]
\newtheorem{lemma}[theorem]{Lemma}
\newtheorem{proposition}[theorem]{Proposition}
\newtheorem{corollary}[theorem]{Corollary}
\theoremstyle{definition}

\theoremstyle{remark}
\newtheorem{remark}[theorem]{Remark}

\numberwithin{equation}{section}

\newcommand{\R}{\mathfrak R}
\newcommand{\I}{\mathfrak I}
\newcommand{\Lh}{\mathscr L}
\newcommand{\Q}{\mathbb Q}
\newcommand{\Z}{\mathbb Z}

\newcommand{\F}{\mathbb F}
\newcommand{\OK}{\mathcal O_K}
\newcommand{\leg}[2]{\left(\frac{#1}{#2}\right)}

\begin{document}

\title[The equation $Cx^2+p^{2m}=y^q$]{On the Diophantine equation $Cx^2+p^{2m}=y^q$}
\author[U. K. Sangale]{Usha K. Sangale}
\address{Swami Ramanand Teerth Marathwada University, Nanded, Maharashtra 431606, India}
\email{ushas073@gmail.com}

\author[P. R. Satpute]{Priyanka R. Satpute}
\address{Swami Ramanand Teerth Marathwada University, Nanded, Maharashtra 431606, India}
\email{priyanka.satpute567@gmail.com}

\begin{abstract}
We classify the primitive positive solutions of $Cx^2+p^{2m}=y^q$ under the hypotheses that $C>3$ is odd and square-free, $p,q$ are odd primes, $p\nmid C$, $q\nmid h(-C)$, and $y$ is odd. The solutions form six explicit families: one for arbitrary $q$ and five exceptional families with $q=3$, three of which arise from half-integral elements when $C\equiv3\pmod8$. The proof combines ideal factorization with Lehmer sequences and shows that representations with second coordinate $v=p^r$, $r>0$, contribute only when $q=3$.
\end{abstract}

\maketitle

\section{Introduction}

Equations of the form
\begin{equation}\label{gen}
C_1x^2+C_2=y^n,\qquad n\ge3,
\end{equation}
with fixed positive coefficients $C_1,C_2$ and unknowns $x,y,n$, are commonly called generalized Lebesgue--Ramanujan--Nagell equations. They include the classical equations of Lebesgue and of Ramanujan and Nagell, as well as many families in which one coefficient is a prime power; see the survey of Le and Soydan \cite{LS} for background and references, and Bugeaud, Mignotte and Siksek \cite{BMS} for the modular approach to the case $C_1=1$.

Several results are particularly close to the problem considered here. Arif and Abu Muriefah \cite{AAM} studied $x^2+q^{2k}=y^n$, and B\'erczes and Pink \cite{BP} determined primitive solutions in the corresponding family for primes $q<100$. Abu Muriefah \cite{AM08} considered
\[
P X^2+Q^{2m}=Y^P,
\]
where $P$ and $Q$ are odd primes. The paper claims a general nonexistence conclusion in part of the range $P>3$, $P\not\equiv7\pmod8$, but the proof of that conclusion has a gap. As Goss \cite{Goss} noted, and Laradji, Mignotte and Tzanakis \cite[\S1]{LMT} explain explicitly, the factorization argument only yields an intermediate coefficient of the form $\pm Q^j$ for some $0\le j\le m$; it does not justify replacing this by $\pm Q^m$. We do not use the unsupported nonexistence conclusion. For positive solutions with $\gcd(X,Y)=1$, Laradji, Mignotte and Tzanakis derive necessary parametric identities when $P>3$ and $P\not\equiv7\pmod8$ (their Proposition 2.1), prove nonexistence when $P\equiv3\pmod8$ (their Theorem 2.3), and give a more detailed analysis in the special case $P=5$. Chakraborty, Hoque and Srinivas \cite{CHS} classified the solutions of the related equation
\begin{equation}\label{CHSeq}
Cx^2+p^{2m}=4y^n,
\end{equation}
under the hypotheses
\[
C\text{ square-free},\quad p\nmid C,\quad m\ge0,\quad n\ge3,
\quad \gcd(x,y)=1,\quad \gcd(n,2h(-C))=1.
\]

Their result is not simply the present theorem with a factor $4$ removed. The factor $4$ changes the parity constraints and the relevant integral representation. Indeed, when $C$ and $p$ are odd, reduction of \eqref{CHSeq} modulo $4$ forces
\[
Cx^2\equiv -p^{2m}\equiv -1\equiv3\pmod4,
\]
so $x$ is odd and $C\equiv3\pmod4$. By contrast, in our equation, the assumptions that $C,p,y$ are odd and that $y$ is odd imply that $x$ is even; the equation then gives $y\equiv1\pmod4$. Thus the two equations lead to different parity regimes and different representations in the relevant quadratic order.

The results also have different scopes. The theorem of \cite{CHS} allows composite exponents $n$ and the case $m=0$, whereas we consider an odd prime exponent $q$ and assume $m\ge1$. Our contribution is an explicit six-family classification for this distinct equation and parity regime, including the exceptional half-integral representations. For prime $q$, the class-number condition $\gcd(q,2h(-C))=1$ in \cite{CHS} is equivalent to $q\nmid h(-C)$; this shared condition, however, does not make the two equations interchangeable.

Patel \cite{Patel} developed a Lucas--Lehmer approach to \eqref{gen}. Writing $C_1C_2=cd^2$ with $c$ square-free and assuming $C_1C_2\not\equiv7\pmod8$ (which forces $y$ odd), she proved that for a prime exponent $n$ one has $n\le5$, or $n=7$ and $y\in\{3,5,9\}$, or $n\mid h(-c)$, or $n\mid\ell-\leg{-c}{\ell}$ for some prime $\ell\mid d$ with $\ell\nmid2c$. She also explained how to solve the equation for any fixed exponent.

We classify the primitive positive solutions of
\[
Cx^2+p^{2m}=y^q,\qquad \gcd(x,y)=1,\quad y\ \text{odd},
\]
where $C>3$ is odd and square-free, $p,q$ are odd primes, $p\nmid C$, and $q\nmid h(-C)$; here $h(-C)$ denotes the class number of $K=\mathbb Q(\sqrt{-C})$. When $C=q$, this specializes, after interchanging the names of the primes, to the equation studied in \cite{AM08,LMT}. It is also the analogue of \eqref{CHSeq} without the factor $4$.

For integers $u,v$ and an odd prime $q$, put
\begin{align}
\R(C,u,v,q)
 &=\sum_{j=0}^{(q-1)/2}\binom{q}{2j}
 u^{q-2j-1}C^{(q-1)/2-j}(-v^2)^j, \label{Rdef}\\
\I(C,u,v,q)
 &=\sum_{j=0}^{(q-1)/2}\binom{q}{2j+1}
 u^{q-2j-1}C^{(q-1)/2-j}(-v^2)^j. \label{Idef}
\end{align}
They allow us to state the classification compactly.

\begin{theorem}\label{mainthm}
Let $C>3$ be an odd square-free integer, and let $p,q$ be odd primes with $p\nmid C$ and $q\nmid h(-C)$. Suppose that
\begin{equation}\label{eqthm}
Cx^2+p^{2m}=y^q,\qquad x,y,m\ge1,\qquad \gcd(x,y)=1,\qquad y\ \text{odd}.
\end{equation}
Then $x$ is even, $y\equiv1\pmod 4$, and there is a unique positive integer $u$ such that exactly one of the following six statements holds.

\smallskip
\noindent\textup{(I)} $u$ is even and
\[
p^m=\bigl|\I(C,u,1,q)\bigr|,\qquad x=u\bigl|\R(C,u,1,q)\bigr|,\qquad y=Cu^2+1 .
\]
\noindent\textup{(II)} $q=3$, $p\ne3$, $u$ is even and
\[
p^{2m}=3Cu^2+1,\qquad x=u\bigl(8Cu^2+3\bigr),\qquad y=4Cu^2+1 .
\]
\noindent\textup{(III)} $q=p=3$, $m\ge2$, $u$ is even and
\[
Cu^2=3^{2m-3}+1,\qquad x=u\bigl(8Cu^2-9\bigr),\qquad y=4Cu^2-3 .
\]
\noindent\textup{(I$'$)} $q=3$, $C\equiv3\pmod8$, $u$ is odd and
\[
p^m=\frac{3Cu^2-1}{8},\qquad x=\frac{u\,(Cu^2-3)}{8},\qquad y=\frac{Cu^2+1}{4} .
\]
\noindent\textup{(II$'$)} $q=3$, $C\equiv3\pmod8$, $p\ne3$, $u$ is odd and
\[
p^{2m}=3Cu^2-8,\qquad x=u\bigl(Cu^2-3\bigr),\qquad y=Cu^2-2 .
\]
\noindent\textup{(III$'$)} $q=p=3$, $C\equiv3\pmod8$, $m\ge2$, $u$ is odd, and for some $\varepsilon\in\{\pm1\}$
\[
Cu^2=3^{2m-3}+8\varepsilon,\qquad x=u\bigl(Cu^2-9\varepsilon\bigr),\qquad y=Cu^2-6\varepsilon .
\]

Conversely, let $C>3$ be odd and square-free, $p$ and $q$ odd primes and $u,m\ge1$. If the conditions on $(q,p,m,u)$ and the displayed equation for $p^m$, $p^{2m}$ or $Cu^2$ in one of \textup{(I)}--\textup{(III$'$)} hold, then the values of $x$ and $y$ given there are positive integers that satisfy \eqref{eqthm}, and $p\nmid C$. The converse does not use the hypothesis $q\nmid h(-C)$.
\end{theorem}

\begin{remark}\label{redundant}
Several side conditions in Theorem~\ref{mainthm} are consequences of the displayed equations and are listed only for clarity. Since $p$ is odd, $u$ is even in (I) (otherwise $\I(C,u,1,q)\equiv\sum_j\binom{q}{2j+1}=2^{q-1}\equiv0\pmod2$) and in (II); in (III), $Cu^2\equiv4\pmod8$. In (I$'$)--(III$'$) the displayed equation forces $Cu^2$ to be odd and $Cu^2\equiv3\pmod8$, hence $u$ odd and $C\equiv3\pmod8$. Finally, $p\ne3$ in (II), (I$'$) and (II$'$), since $3Cu^2+1$, $3Cu^2-1$ and $3Cu^2-8$ are prime to $3$.
\end{remark}

Theorem~\ref{mainthm} gives a complete parametrization rather than only a bound on the exponent or an algorithm for a fixed exponent. Patel's result \cite{Patel} provides broad restrictions for the generalized equation $C_1x^2+C_2=y^n$ and an effective procedure for fixed data; the present theorem gives closed-form families for the specific coefficient pattern $C_1=C$ and $C_2=p^{2m}$. In the special case $C=q$, Corollary~\ref{LMTcor} gives the same necessary polynomial parametrization as Laradji, Mignotte and Tzanakis \cite[Proposition~2.1]{LMT} in their range $q>3$, $q\not\equiv7\pmod8$, after interchanging the names of the primes. Their proposition is a parametrization that every solution must satisfy, not a complete resolution for all exponents. The present corollary also covers $q\equiv7\pmod8$ provided $y$ is odd, a case excluded from their proposition. Their Theorem~2.3 separately proves nonexistence when the exponent is $q\equiv3\pmod8$ and $q>3$. The comparison with Chakraborty, Hoque and Srinivas \cite{CHS} is also relevant: their equation has a factor $4$ on the right, whereas the present equation does not.

The proof identifies the role of the second coordinate $v$ in the representation supplied by the ideal-factorization lemma. Since $v=p^r$ for some $0\le r\le m$, the case $r=0$ gives the generic family, while $r>0$ is ruled out for $q\ge5$ and leads to the additional families when $q=3$. In the special case $C=q$, this keeps track of the intermediate exponent that was not justified in the original argument discussed in \cite{AM08,LMT}. The three further families (I$'$)--(III$'$) arise from half-integral elements of $\mathcal O_K$ and occur only when $q=3$ and $C\equiv3\pmod8$; they would be missed by working solely in $\mathbb Z[\sqrt{-C}]$.

A further consequence is the congruence $p\equiv\leg{-C}{p}\pmod q$ for the generic families (Proposition~\ref{congruence}), analogous to a condition in Patel \cite[Theorem~1]{Patel}. The hypothesis that $y$ is odd is essential only when $C\equiv7\pmod8$; Remark~\ref{yeven} gives examples in that case. The ideal-factorization step itself is standard; Lemma~\ref{replemma} adapts it to the maximal order and makes the half-integral case explicit (Remark~\ref{yuanremark}).

The remainder of the paper is organized as follows. Section~\ref{sec:prelim} contains the representation lemma and the auxiliary results, Section~\ref{sec:main} the proof of Theorem~\ref{mainthm}, Section~\ref{sec:cons} consequences, and Section~\ref{sec:examples} examples and a computational check. Section~\ref{sec:final} concludes with further questions.

\section{Preliminaries}\label{sec:prelim}

\subsection{Notation}
Throughout, $C>3$ is an odd square-free integer, $K=\Q(\sqrt{-C})$, $\OK$ is its ring of integers and $h(-C)$ its class number. Since $C>3$, the only roots of unity in $K$ are $\pm1$. The polynomials $\R(C,u,v,q)$ and $\I(C,u,v,q)$ of \eqref{Rdef}--\eqref{Idef} depend only on $u^2$ and $v^2$, and for $\mu\in\{\pm1\}$ the binomial theorem gives
\begin{equation}\label{powerexp}
 (u\sqrt C+\mu vi)^q
 =u\,\R(C,u,v,q)\sqrt C+\mu v\,\I(C,u,v,q)\,i .
\end{equation}
Since $v+u\sqrt{-C}=i\,(u\sqrt C-vi)$ and $i^{q}=(-1)^{(q-1)/2}i$, \eqref{powerexp} with $\mu=-1$ gives
\begin{equation}\label{powerexpK}
(v+u\sqrt{-C})^q=(-1)^{(q-1)/2}\bigl(v\,\I(C,u,v,q)+u\,\R(C,u,v,q)\sqrt{-C}\bigr).
\end{equation}
For $q=3$,
\begin{equation}\label{RI3}
\R(C,u,v,3)=Cu^2-3v^2,\qquad \I(C,u,v,3)=3Cu^2-v^2 .
\end{equation}
The last term of \eqref{Idef} is $(-1)^{(q-1)/2}v^{q-1}$ and every other term is divisible by $Cu^2$; similarly, the last term of \eqref{Rdef} is $(-1)^{(q-1)/2}qv^{q-1}$. Hence
\begin{equation}\label{Icong}
\I(C,u,1,q)\equiv\pm1,\qquad \R(C,u,1,q)\equiv\pm q\pmod{Cu^2}.
\end{equation}

\subsection{The parity of \texorpdfstring{$y$}{y}}

\begin{lemma}\label{yodd}
Let $x,y,m\ge1$ satisfy $Cx^2+p^{2m}=y^q$ with $\gcd(x,y)=1$, where $p,q$ are odd primes. If $y$ is even, then $x$ is odd and $C\equiv7\pmod8$. In particular $y$ is odd whenever $C\not\equiv7\pmod 8$.
\end{lemma}

\begin{proof}
If $y$ is even then $x$ is odd, so $Cx^2+p^{2m}\equiv C+1\pmod 8$, while $8\mid y^q$ because $q\ge3$. Hence $C\equiv7\pmod 8$.
\end{proof}

\subsection{A representation lemma}

\begin{lemma}\label{replemma}
Let $q$ be an odd prime with $q\nmid h(-C)$, and let $x,Y,Z$ be positive integers such that
\[
Cx^2+Y^2=Z^q,\qquad Z\ \text{odd},\qquad \gcd(Y,Z)=1 .
\]
Then there are positive integers $u,v$ with $\gcd(Cu,v)=1$ and signs $\delta,\mu\in\{\pm1\}$ such that one of the following holds:
\begin{enumerate}
\item[(a)] $x\sqrt C+Yi=\delta(u\sqrt C+\mu vi)^q$ and $Z=Cu^2+v^2$;
\item[(b)] $q=3$, $C\equiv3\pmod 8$, $u$ and $v$ are odd,
$x\sqrt C+Yi=\delta\bigl((u\sqrt C+\mu vi)/2\bigr)^3$ and $Z=(Cu^2+v^2)/4$.
\end{enumerate}
Moreover, the choice between \textup{(a)} and \textup{(b)}, as well as the pair $(u,v)$, is uniquely determined by $(x,Y)$.
\end{lemma}

\begin{proof}
In $\OK$ we have the factorization
\begin{equation}\label{ideal-factor}
(Y+x\sqrt{-C})(Y-x\sqrt{-C})=Z^q.
\end{equation}
We first show that the two factors on the left generate coprime ideals. Suppose that a prime ideal $\mathfrak P$ divides both. By \eqref{ideal-factor}, it divides $Z^q$, so it lies above an odd prime $\ell\mid Z$, because $Z$ is odd. It also divides the sum of the two factors, namely $2Y$. Since $2\notin\mathfrak P$, it follows that $Y\in\mathfrak P$, and hence $\ell\mid Y$. This contradicts $\gcd(Y,Z)=1$.

The two principal ideals are therefore coprime. Their product is $(Z)^q$, so unique factorization of ideals gives an ideal $\mathfrak a$ such that
\[
(Y+x\sqrt{-C})=\mathfrak a^q.
\]
The class of $\mathfrak a$ has order dividing $q$. It also has order dividing $h(-C)$, and the hypothesis $q\nmid h(-C)$ therefore forces $\mathfrak a$ to be principal. Write $\mathfrak a=(\gamma)$. Since the only units in $\OK$ are $\pm1$ and $q$ is odd, we may replace $\gamma$ by $-\gamma$ if necessary to obtain
\begin{equation}\label{gammaeq}
Y+x\sqrt{-C}=\gamma^q,
\qquad
\gamma=\frac{s+t\sqrt{-C}}{2},
\qquad s,t\in\Z,\quad s\equiv t\pmod2.
\end{equation}
Taking norms gives
\[
Z=N(\gamma)=\frac{s^2+Ct^2}{4}.
\]
Moreover, neither $s$ nor $t$ is zero. If $t=0$, then $\gamma^q\in\Z$, which would force $x=0$; if $s=0$, then $\gamma^q\in\Z\sqrt{-C}$, which would force $Y=0$. Both contradict the hypotheses.

Suppose that $s$ and $t$ are odd. Then $\gamma\notin\Z[\sqrt{-C}]$, so $C\equiv3\pmod4$. If $C\equiv7\pmod8$, then
\[
s^2+Ct^2\equiv1+7\equiv0\pmod8,
\]
which would make $Z=(s^2+Ct^2)/4$ even. This contradicts the hypothesis that $Z$ is odd. Consequently,
\[
C\equiv3\pmod8.
\]
In this case $2$ is inert in $K$, and $\OK/2\OK\cong\F_4$. Put
\[
\omega=\frac{1+\sqrt{-C}}2.
\]
Since $\sqrt{-C}=2\omega-1$, the image of $\Z[\sqrt{-C}]$ in $\F_4$ is the prime field $\F_2$. On the other hand,
\[
\gamma=\frac{s-t}{2}+t\omega,
\]
and $t$ is odd. Thus the image $\bar\gamma$ of $\gamma$ is not in $\F_2$. In particular, $\bar\gamma\ne0,1$, so $\bar\gamma$ has order $3$ in $\F_4^\times$.

But $\gamma^q=Y+x\sqrt{-C}\in\Z[\sqrt{-C}]$, so $\bar\gamma^q\in\F_2^\times=\{1\}$. It follows that $3\mid q$, and hence $q=3$. A similar argument appears in the proof of \cite[Proposition~2.1]{LMT}.

Put $u=|t|/2$, $v=|s|/2$ if $s,t$ are even, and $u=|t|$, $v=|s|$ if $s,t$ are odd. Conjugating \eqref{gammaeq} and multiplying by $i$, and using $i(Y-x\sqrt{-C})=x\sqrt C+Yi$, $i(s-t\sqrt{-C})=t\sqrt C+si$ and $i=(-1)^{(q-1)/2}i^q$, we obtain
\[
x\sqrt C+Yi=(-1)^{(q-1)/2}\Bigl(\frac{t\sqrt C+si}{2}\Bigr)^q .
\]
Since
\[
t\sqrt C+si=\operatorname{sgn}(t)\bigl(|t|\sqrt C+\operatorname{sgn}(st)|s|\,i\bigr)
\]
and $q$ is odd, this is (a) or (b) with $\delta=(-1)^{(q-1)/2}\operatorname{sgn}(t)$ and $\mu=\operatorname{sgn}(st)$; the formula for $Z$ is the norm computed above.

To see that $\gcd(Cu,v)=1$, let $\ell$ be a prime dividing $Cu$ and $v$. Comparing coefficients of $i$ by \eqref{powerexp} gives $Y=\pm v\,\I(C,u,v,q)$ in case (a) and $8Y=\pm v\,\I(C,u,v,3)$ in case (b); so $\ell\mid Y$ in case (a), and also in case (b) because $v$ is odd. Moreover $\ell$ divides $Cu^2+v^2$, which is $Z$ in case (a) and $4Z$ in case (b), so $\ell\mid Z$ (in case (b) since $\ell$ is odd). This contradicts $\gcd(Y,Z)=1$.

For uniqueness, suppose that (a) or (b) holds with some $u,v,\delta,\mu$. Replacing $v$ by $\mu v$ in \eqref{powerexpK} and comparing with the coefficients of $x\sqrt C+Yi$ shows that $\gamma'=\pm(\mu v+u\sqrt{-C})$ in case (a), and $\gamma'=\pm(\mu v+u\sqrt{-C})/2$ in case (b), satisfies $\gamma'^q=Y+x\sqrt{-C}$ for a suitable sign. Then $\gamma'/\gamma$ is a root of unity in $K$ whose $q$th power is $1$, so $\gamma'=\gamma$ because $q$ is odd. Hence the case is determined by whether $\gamma\in\Z[\sqrt{-C}]$, and $(u,v)$ is determined by $\gamma$.
\end{proof}

\begin{remark}\label{yuanremark}
Lemma~\ref{replemma} is a version of a standard ideal-factorization argument, adapted to the maximal order, the hypothesis $q\nmid h(-C)$ and odd $Z$; compare Yuan \cite{Yuan}, who treats the more general equation $ax^2+by^2=ck^n$ by means of binary quadratic forms. The hypothesis that $Z$ is odd is used to make $Y\pm x\sqrt{-C}$ coprime. Case (b) cannot be omitted: for $C=11$,
\[
11\cdot2^2+3^4=5^3,\qquad
3^2+2\sqrt{-11}=\left(\frac{-3+\sqrt{-11}}2\right)^3,
\]
and $3^2+2\sqrt{-11}$ is not a cube in $\Z[\sqrt{-11}]$. This case is invisible to an argument in the order $\Z[\sqrt{-C}]$ of discriminant $-4C$, whose class number is $h(-4C)$. For $C\equiv3\pmod4$ the class number formula for orders gives
\[
h(-4C)=\left(2-\leg{-C}{2}\right)h(-C),
\]
where $\leg{-C}{2}$ is the Kronecker symbol; in particular $h(-4C)=3h(-C)$ when $C\equiv3\pmod8$. So for $q=3$ and $C\equiv3\pmod8$ the condition $q\nmid h(-4C)$ never holds, and the half-integral case must be handled explicitly.
\end{remark}

\subsection{Lehmer pairs}
A pair $(\alpha,\beta)$ of algebraic integers is a \emph{Lehmer pair} if $(\alpha+\beta)^2$ and $\alpha\beta$ are non-zero coprime rational integers and $\alpha/\beta$ is not a root of unity. For odd $n$ its $n$th \emph{Lehmer number} is
\[
\Lh_n(\alpha,\beta)=\frac{\alpha^n-\beta^n}{\alpha-\beta},
\]
and for even $n$ one divides by $\alpha^2-\beta^2$ instead. A prime $\ell$ is a \emph{primitive divisor} of $\Lh_n(\alpha,\beta)$ if $\ell\mid\Lh_n(\alpha,\beta)$ but $\ell\nmid(\alpha^2-\beta^2)^2\Lh_1(\alpha,\beta)\cdots\Lh_{n-1}(\alpha,\beta)$. A Lehmer pair is \emph{$n$-defective} if $\Lh_n(\alpha,\beta)$ has no primitive divisor. Two Lehmer pairs $(\alpha_1,\beta_1)$, $(\alpha_2,\beta_2)$ are \emph{equivalent} if $\alpha_1/\alpha_2=\beta_1/\beta_2\in\{\pm1,\pm\sqrt{-1}\}$. Every Lehmer pair is equivalent to $\bigl((\sqrt a+\sqrt b)/2,(\sqrt a-\sqrt b)/2\bigr)$ with $(a,b)=\bigl((\alpha+\beta)^2,(\alpha-\beta)^2\bigr)$, and passing to an equivalent pair changes $(a,b)$ at most to $(-a,-b)$.

\begin{theorem}[Bilu--Hanrot--Voutier \cite{BHV}; Voutier \cite{Voutier95}]\label{bhv}
If $n>30$, there are no $n$-defective Lehmer pairs. If $q$ is a prime with $7\le q\le 29$ and $(\alpha,\beta)$ is $q$-defective, then $q\in\{7,13\}$ and, up to equivalence, $(a,b)$ is one of
\[
(1,-7),\ (1,-19),\ (3,-5),\ (5,-7),\ (13,-3),\ (14,-22)\quad(q=7);\qquad (1,-7)\quad(q=13).
\]
\end{theorem}

This is \cite[Theorem~1.4]{BHV} together with Table~2 of \cite{BHV} (see also \cite{Voutier95}). The recent corrections of Voutier \cite{VoutierCorr} concern only the indices $n\in\{3,4,5,6,8,10,12\}$; since we treat $q=3$ and $q=5$ directly, they do not affect our arguments.

\subsection{Lucas numbers}
Let $(F_k)$ and $(L_k)$ denote the Fibonacci and Lucas sequences. The integer solutions of $X^2-5Y^2=\pm4$ are exactly $(X,Y)=(\pm L_k,\pm F_k)$, $k\ge0$, and the only Lucas numbers that are perfect squares are $L_1=1$ and $L_3=4$ (Cohn \cite{Cohn}).

\section{Proof of the main theorem}\label{sec:main}

\subsection*{Step 1: parity and coprimality}
We first show that $p$ does not divide $y$. Indeed, if $p\mid y$, then the Diophantine equation gives $p\mid Cx^2$. Since $p\nmid C$, it follows that $p\mid x$, contrary to $\gcd(x,y)=1$. Thus
\begin{equation}\label{coprime}
p\nmid y.
\end{equation}

Next, $C$, $p$ and $y$ are odd. Consequently, $Cx^2+p^{2m}=y^q$ implies that $Cx^2$ is even, so $x$ is even. In particular, $4\mid Cx^2$. Reducing the equation modulo $4$ now gives
\[
y^q\equiv p^{2m}\equiv1\pmod4.
\]
Since $q$ is odd, this yields $y\equiv1\pmod4$.

\subsection*{Step 2: parametrization}
By \eqref{coprime}, Lemma~\ref{replemma} applies with $Y=p^m$ and $Z=y$; let $u,v,\delta,\mu$ be the resulting parameters. Families (I)--(III) will arise from case (a) of Lemma~\ref{replemma} and families (I$'$)--(III$'$) from case (b). Comparing coefficients by \eqref{powerexp} and \eqref{RI3},
\begin{align}
x&=\delta u\,\R(C,u,v,q),& p^m&=\delta\mu v\,\I(C,u,v,q) &&\text{in case (a)},\label{realimag}\\
8x&=\delta u\,(Cu^2-3v^2),& 8p^m&=\delta\mu v\,(3Cu^2-v^2) &&\text{in case (b)}.\label{halfrealimag}
\end{align}
In case (a), $v\mid p^m$; in case (b), $v$ is odd and $v\mid 8p^m$. In both cases
\[
v=p^r,\qquad 0\le r\le m .
\]
Thus $r=0$, $r=m$, and $0<r<m$ correspond to the generic family, the second family, and the third family, respectively, in each of the integral and half-integral cases.

\subsection*{Case (a)}
Since $v$ and $y=Cu^2+v^2$ are odd, $u$ is even. If $r=0$, then $v=1$ and \eqref{realimag} is family (I). Assume from now on that $r\ge1$.

\emph{The Lehmer pair.} Put
\[
\alpha=u\sqrt C+\mu vi,\qquad \beta=u\sqrt C-\mu vi.
\]
Then
\begin{equation}\label{lehmerdata}
(\alpha+\beta)^2=4Cu^2,\qquad
\alpha\beta=Cu^2+v^2=y,\qquad
(\alpha-\beta)^2=-4v^2.
\end{equation}
Since $y$ is odd and $\gcd(Cu,v)=1$, we have $\gcd(4Cu^2,y)=1$. Also,
\[
\frac{\alpha}{\beta}=\frac{\alpha^2}{y}
=\frac{Cu^2-v^2+2\mu uv\sqrt{-C}}{y}\in K.
\]
This quotient is not $\pm1$, because $uv\ne0$. The only roots of unity in $K$ are $\pm1$, so $\alpha/\beta$ is not a root of unity. Thus $(\alpha,\beta)$ is a Lehmer pair.

The representation from Lemma~\ref{replemma} gives
\[
x\sqrt C+p^mi=\delta\alpha^q,
\qquad
x\sqrt C-p^mi=\delta\beta^q.
\]
Subtracting these identities and using $\alpha-\beta=2\mu vi$ yields
\[
\alpha^q-\beta^q=2\delta p^mi,
\qquad
\alpha-\beta=2\mu vi.
\]
Therefore
\begin{equation}\label{Lq}
\Lh_q(\alpha,\beta)=\delta\mu\,p^{m-r}.
\end{equation}
Every prime divisor of this Lehmer number is $p$. Since $r\ge1$, we also have $p\mid v$, and hence
\[
p\mid(\alpha^2-\beta^2)^2=-16Cu^2v^2.
\]
It follows that $(\alpha,\beta)$ is $q$-defective.

\emph{The case $q\ge7$.} By Theorem~\ref{bhv}, $q\in\{7,13\}$ and $|a|=4Cu^2$ belongs to $\{1,3,5,13,14\}$. This is impossible, since $u\ge2$ and $C\ge5$ give $4Cu^2\ge80$.

\emph{The case $q=5$.} With $A=(\alpha+\beta)^2=4Cu^2$ and $B=\alpha\beta=y$ we have $\Lh_5(\alpha,\beta)=\alpha^4+\alpha^3\beta+\alpha^2\beta^2+\alpha\beta^3+\beta^4=A^2-3AB+B^2$. With $w=Cu^2$ and $B=w+v^2$ this becomes
\begin{equation}\label{L5}
\Lh_5(\alpha,\beta)=5w^2-10wv^2+v^4=5X^2-4v^4,\qquad X=w-v^2 .
\end{equation}
Note that $X$ is odd. Suppose first that $r<m$. By \eqref{Lq}, the prime $p$ divides $\Lh_5(\alpha,\beta)$, and by assumption $p\mid v$. Reducing the expression in \eqref{L5} modulo $p$ gives $p\mid5X^2$. Since $p\nmid Cu$ and $X=Cu^2-v^2\equiv Cu^2\not\equiv0\pmod p$, we obtain $p=5$.

Because $5\mid v$ but $5\nmid w$, the expression
\[
\Lh_5(\alpha,\beta)=5w^2-10wv^2+v^4
\]
is divisible by $5$ but not by $25$. Equation~\eqref{Lq} therefore implies $m-r=1$, and \eqref{L5} becomes
\[
5X^2-4v^4=\pm5.
\]
As $v=5^r$, the two possible signs give, respectively,
\[
X^2-1=4\cdot5^{4r-1}
\qquad\text{or}\qquad
X^2+1=4\cdot5^{4r-1}.
\]
The second equation is impossible modulo $4$. In the first equation,
\[
\frac{X-1}{2}\,\frac{X+1}{2}=5^{4r-1}.
\]
The two factors are coprime consecutive integers. Their product is a positive power of $5$, so one factor must be $\pm1$ and the other $\pm5^{4r-1}$; this is incompatible with their being consecutive.

It remains to consider $r=m$. In this case, \eqref{L5} and \eqref{Lq} give
\[
5X^2-4v^4=\pm1,
\]
or equivalently,
\[
(4v^2)^2-5(2X)^2=\mp4.
\]
The description of the solutions of $U^2-5V^2=\pm4$ implies that $4v^2$ is a Lucas number. By Cohn's theorem, the only Lucas numbers that are squares are $1$ and $4$. Since $4v^2$ is a square divisible by $4$, we must have $4v^2=4$, so $v=1$, contradicting $r\ge1$. Thus $q=5$ is impossible.

\emph{The case $q=3$.} Here $\Lh_3(\alpha,\beta)=(\alpha+\beta)^2-\alpha\beta=3Cu^2-v^2$, so by \eqref{Lq}
\begin{equation}\label{L3}
\bigl|3Cu^2-p^{2r}\bigr|=p^{m-r}.
\end{equation}
Suppose first that $r=m$. Then \eqref{L3} gives
\[
3Cu^2-p^{2m}=\pm1.
\]
Since $4\mid u^2$ and $p^{2m}\equiv1\pmod4$, the left-hand side is congruent to $3$ modulo $4$. Thus the sign is negative, and
\begin{equation}\label{caseIIeq}
p^{2m}=3Cu^2+1.
\end{equation}
In particular, $p\ne3$. Put $A=Cu^2$. Since $v^2=p^{2m}=3A+1$, we obtain
\[
y=A+v^2=4A+1.
\]
Moreover, by \eqref{RI3},
\[
\R(C,u,v,3)=A-3v^2=-(8A+3).
\]
Equation~\eqref{realimag} now gives $x=u(8A+3)$, which is family~(II).

Suppose next that $1\le r<m$. Reducing \eqref{L3} modulo $p$ gives $p\mid3Cu^2$. As $p\nmid Cu$, we conclude that $p=3$. Dividing \eqref{L3} by $3$ yields
\[
\bigl|Cu^2-3^{2r-1}\bigr|=3^{m-r-1}.
\]
The left-hand side is not divisible by $3$, because $3\nmid Cu$. Therefore $m=r+1$, and
\[
Cu^2=3^{2r-1}\pm1.
\]
Now $3^{2r-1}\equiv3\pmod8$, while $4\mid Cu^2$. Only the plus sign is possible, so
\begin{equation}\label{caseIIIeq}
Cu^2=3^{2m-3}+1,\qquad m\ge2.
\end{equation}
Writing $A=Cu^2$, we have $v^2=3^{2r}=3(A-1)$, whence
\[
y=A+v^2=4A-3,
\qquad
\R(C,u,v,3)=A-9(A-1)=-(8A-9).
\]
Thus $x=u(8A-9)$, giving family~(III).

\subsection*{Case (b)}
Now $q=3$, $C\equiv3\pmod8$ (so $C\ge11$), $u,v$ are odd and $y=(Cu^2+v^2)/4$. By \eqref{halfrealimag},
\begin{equation}\label{halfL3}
\bigl|3Cu^2-p^{2r}\bigr|=8p^{m-r}.
\end{equation}

If $r=0$, then $v=1$. Since $C\ge11$ and $u$ is odd, we have $3Cu^2-1>0$. Equation~\eqref{halfL3} therefore gives
\[
8p^m=3Cu^2-1.
\]
Substituting $v=1$ into \eqref{halfrealimag} yields
\[
8x=u(Cu^2-3),\qquad 4y=Cu^2+1.
\]
These are precisely the formulas in family~(I$'$).

If $r=m$, then \eqref{halfL3} gives
\[
3Cu^2-p^{2m}=\pm8.
\]
We must have $p\ne3$, since otherwise the left-hand side is divisible by $3$ but the right-hand side is not. Thus $p^{2m}\equiv1\pmod3$. Since $3Cu^2-p^{2m}\equiv-1\equiv2\pmod3$, the sign must be positive, and hence
\[
p^{2m}=3Cu^2-8.
\]
Put $A=Cu^2$. Then $v^2=3A-8$, so
\[
y=\frac{A+v^2}{4}=A-2.
\]
Also, \eqref{halfrealimag} gives
\[
8x=u\lvert A-3v^2\rvert=8u(A-3),
\]
which yields family~(II$'$).

Finally, suppose that $1\le r<m$. Reducing \eqref{halfL3} modulo $p$ and using $p\nmid Cu$ gives $p=3$. After division by $3$, we obtain
\[
\bigl|Cu^2-3^{2r-1}\bigr|=8\cdot3^{m-r-1}.
\]
Since $3\nmid Cu$, the left-hand side is not divisible by $3$. Hence $m=r+1$, and for some $\varepsilon\in\{\pm1\}$,
\begin{equation}\label{caseIIIprimeeq}
Cu^2=3^{2m-3}+8\varepsilon.
\end{equation}
Writing $A=Cu^2$, we have
\[
v^2=3^{2r}=3(A-8\varepsilon),
\qquad
y=\frac{A+v^2}{4}=A-6\varepsilon.
\]
Furthermore,
\[
8x=u\lvert A-3v^2\rvert
=u\lvert A-9(A-8\varepsilon)\rvert
=8u(A-9\varepsilon).
\]
The last equality follows from $A\ge11$, which ensures $A-9\varepsilon>0$. This is family~(III$'$).

\subsection*{Converse}
Let the data satisfy the conditions of one of the families and put $A=Cu^2$. By Remark~\ref{redundant}, $A\equiv3\pmod8$ in (I$'$)--(III$'$); hence in (I$'$) the numbers $x=u(A-3)/8$ and $y=(A+1)/4$ are integers. For (I), \eqref{powerexp} gives
\[
C\bigl(u\R(C,u,1,q)\bigr)^2+\I(C,u,1,q)^2=\bigl|(u\sqrt C+i)^q\bigr|^2=(Cu^2+1)^q .
\]
For the remaining families the following polynomial identities in $A$ give the equation:
\begin{align*}
&\text{(II)}\quad A(8A+3)^2+(3A+1)=(4A+1)^3,\\
&\text{(III)}\quad A(8A-9)^2+27(A-1)=(4A-3)^3,\qquad 3^{2m}=27\cdot3^{2m-3}=27(A-1),\\
&\text{(I$'$)}\quad A(A-3)^2+(3A-1)^2=(A+1)^3\quad\text{(the equation multiplied by $64$)},\\
&\text{(II$'$)}\quad A(A-3)^2+(3A-8)=(A-2)^3,\\
&\text{(III$'$)}\quad A(A-9\varepsilon)^2+27(A-8\varepsilon)=(A-6\varepsilon)^3,\qquad 3^{2m}=27(A-8\varepsilon).
\end{align*}

\emph{Positivity.} In (I), $\R(C,u,1,q)\equiv\pm q\pmod{Cu^2}$ by \eqref{Icong}; since $Cu^2>1$ is not prime ($4\mid u^2$), it does not divide $q$, so $\R(C,u,1,q)\ne0$ and $x>0$. In (II) $x>0$ trivially; in (III), $A\ge5\cdot2^2$, so $8A-9>0$; in (I$'$)--(III$'$), $A\ge11$, so $A-3>0$ and $A-9\varepsilon>0$.

\emph{Parity of $y$ and $p\nmid C$.} In (I)--(III), $y$ is odd because $u$ is even; in (I$'$)--(III$'$) because $A\equiv3\pmod8$. In (I), $p^m\equiv\pm1\pmod C$ by \eqref{Icong}; in (II), (I$'$) and (II$'$), $p^{2m}$ or $8p^m$ is congruent to $1$, $-1$ or $-8$ modulo the odd number $C$; in (III) and (III$'$), $p=3$ and $3\nmid A$. In each case $p\nmid C$.

\emph{Coprimality.} If a prime $\ell$ divides $x$ and $y$, then $\ell\mid p^{2m}$, so $\ell=p$; hence it suffices to show that $p\nmid y$. In (II), $3y=4p^{2m}-1$; in (II$'$), $3y=p^{2m}+2$; in (I$'$), $12y=8p^m+4$; in (III) and (III$'$), $p=3$ and $y\equiv A\not\equiv0\pmod 3$. In each case $p\nmid y$. In (I), suppose that $p\mid y$, so that $p\mid x$, and put $\zeta=1+u\sqrt{-C}\in\OK$. By \eqref{powerexpK}, $\zeta^q=\pm\bigl(\I(C,u,1,q)+u\R(C,u,1,q)\sqrt{-C}\bigr)$, so $p^m+x\sqrt{-C}$ is one of $\pm\zeta^q$, $\pm\bar\zeta^q$, and $p$ divides both $\zeta^q$ and $\bar\zeta^q$ in $\OK$. A prime ideal above $p$ then divides $\zeta$ and $\bar\zeta$, hence $\zeta+\bar\zeta=2$, which is absurd.

None of these verifications uses $q\nmid h(-C)$.

\subsection*{Uniqueness and exclusivity}
For a solution in one of the six families put $v=1$ in (I) and (I$'$), $v=p^m$ in (II) and (II$'$), and $v=3^{m-1}$ in (III) and (III$'$). The computations above show that $x=|u\,\R(C,u,v,q)|$ and $p^m=|v\,\I(C,u,v,q)|$ in (I)--(III), and $8x=|u\,\R(C,u,v,3)|$ and $8p^m=|v\,\I(C,u,v,3)|$ in (I$'$)--(III$'$); so, by \eqref{powerexp} and a suitable choice of signs, $x\sqrt C+p^mi$ is represented as in case (a) of Lemma~\ref{replemma} for (I)--(III) and as in case (b) for (I$'$)--(III$'$). By the uniqueness part of Lemma~\ref{replemma}, the case and the pair $(u,v)$ are determined by the solution. The six families are distinguished by the case together with $v$ (and in (III$'$) the sign $\varepsilon$ is determined by $Cu^2$), so they are pairwise exclusive and $u$ is unique. This completes the proof of Theorem~\ref{mainthm}. \qed

\section{Consequences}\label{sec:cons}

\begin{corollary}\label{corgeneric}
Under the hypotheses of Theorem~\ref{mainthm}, if $q\ge5$ then every solution of \eqref{eqthm} belongs to family \textup{(I)}, and if $q=3$ and $C\not\equiv3\pmod 8$ then every solution belongs to one of the families \textup{(I)}, \textup{(II)}, \textup{(III)}.
\end{corollary}

\begin{proposition}\label{congruence}
Under the hypotheses of Theorem~\ref{mainthm}, if a solution belongs to family \textup{(I)} or \textup{(I$'$)}, then
\[
p\equiv\leg{-C}{p}\pmod q ,
\]
where $\leg{\cdot}{p}$ is the Legendre symbol. In particular $p\equiv\pm1\pmod q$.
\end{proposition}

\begin{proof}
We first note that $p\ne q$. In family (I), every binomial coefficient in \eqref{Idef} except the last is divisible by $q$, so $\I(C,u,1,q)\equiv\pm1\pmod q$ and $q\nmid p^m$. In family (I$'$), $8p^m=3Cu^2-1\not\equiv0\pmod3$.

Let $\alpha=u\sqrt C+\mu i$ in family (I) and $\alpha=(u\sqrt C+\mu i)/2$ in family (I$'$), and let $\beta=\bar\alpha$. In both cases $\alpha\beta=y$ and, as in the proof of Theorem~\ref{mainthm}, $\alpha^q-\beta^q=(\alpha-\beta)\cdot c\,p^m$ with $c\in\{\pm1,\pm2\}$; so $\alpha^q-\beta^q$ is divisible by $p$ in the ring of integers of $K(i)$. Put $\gamma=\alpha/\beta=\alpha^2/y\in K$; then $N_{K/\Q}(\gamma)=\gamma\bar\gamma=1$. Let $\mathfrak p$ be a prime of $\OK$ above $p$. Since $\alpha\beta=y$ and $p\nmid y$ by \eqref{coprime}, $\alpha$ and $\beta$ are units at every prime of $K(i)$ above $\mathfrak p$. Hence $\gamma$ is a unit at $\mathfrak p$, and $\gamma^q-1=(\alpha^q-\beta^q)/\beta^q$ has positive valuation at the primes of $K(i)$ above $\mathfrak p$; as $\gamma^q-1\in K$, this gives $\gamma^q\equiv1\pmod{\mathfrak p}$.

To see that $\gamma\not\equiv1\pmod{\mathfrak p}$, note that in both families
\[
\gamma-1=\frac{2(\mu u\sqrt{-C}-1)}{Cu^2+1},
\]
where $Cu^2+1$, which equals $y$ in family (I) and $4y$ in family (I$'$), is prime to $p$. If $\gamma\equiv1\pmod{\mathfrak p}$, then $\mu u\sqrt{-C}\equiv1\pmod{\mathfrak p}$, and squaring gives $-Cu^2\equiv1\pmod{\mathfrak p}$, so $p\mid Cu^2+1$, a contradiction. Hence $\gamma$ has order exactly $q$ in $(\OK/\mathfrak p)^\times$.

As $p$ is odd and $p\nmid C$, $p$ is unramified in $K$. If $\leg{-C}{p}=1$, then $\OK/\mathfrak p\cong\F_p$ and $q\mid p-1$. If $\leg{-C}{p}=-1$, then $\OK/\mathfrak p\cong\F_{p^2}$, the Frobenius acts as the nontrivial automorphism of $K$, and $\gamma^{p+1}\equiv\gamma\bar\gamma=1\pmod{\mathfrak p}$, so $q\mid p+1$.
\end{proof}

\begin{remark}
Proposition~\ref{congruence} is the analogue for \eqref{eqthm} of condition (iv) of \cite[Theorem~1]{Patel}, here valid for every $q$ and without the restriction $C\not\equiv7\pmod8$. It gives a quick sieve: if $q\ge5$ and $p\not\equiv\pm1\pmod q$, then \eqref{eqthm} has no solution. For families (II) and (III) the congruence need not hold; for instance $p=3$ in (III).
\end{remark}

\begin{corollary}\label{LMTcor}
Let $q>3$ be a prime and $p\ne q$ an odd prime. If $qx^2+p^{2m}=y^q$ with $x,y,m\ge1$, $\gcd(x,y)=1$ and $y$ odd, then there is an even $u\ge2$ with
\[
p^m=\bigl|\I(q,u,1,q)\bigr|,\qquad x=u\bigl|\R(q,u,1,q)\bigr|,\qquad y=qu^2+1,
\]
and $p\equiv\leg{-q}{p}\pmod q$.
\end{corollary}

\begin{proof}
Take $C=q$. Then $h(-q)<q$ (see the proof of \cite[Proposition~2.1]{LMT}), so $q\nmid h(-q)$ and Theorem~\ref{mainthm} applies; since $q>3$, only family (I) is possible, and Proposition~\ref{congruence} gives the congruence.
\end{proof}

When $q\not\equiv7\pmod8$, Lemma~\ref{yodd} makes the parity condition on $y$ automatic, and Corollary~\ref{LMTcor} agrees with the necessary parametrization in \cite[Proposition~2.1]{LMT}; for $q=5$ and $u=2$, it gives $5\cdot410^2+1801^2=21^5$. If $q\equiv3\pmod8$ and $q>3$, the nonexistence theorem of Laradji, Mignotte and Tzanakis \cite[Theorem~2.3]{LMT} shows that no solution occurs, so the corresponding parametric conditions cannot be realized. Corollary~\ref{LMTcor} also covers $q\equiv7\pmod8$ when $y$ is odd.

\begin{remark}\label{yeven}
The hypothesis that $y$ is odd cannot be omitted in Theorem~\ref{mainthm}. By Lemma~\ref{yodd} it only matters when $C\equiv7\pmod 8$, and then there are solutions with $y$ even, $x$ odd, $p\nmid C$ and $q\nmid h(-C)$, for example
\[
7\cdot1^2+5^2=2^5,\quad 7\cdot1^2+11^2=2^7,\quad 7\cdot1^2+181^2=32^3,\quad 15\cdot1^2+7^2=4^3,\quad 23\cdot1^2+3^2=2^5,
\]
with $h(-7)=1$, $h(-15)=2$ and $h(-23)=3$; they violate the conclusion that $x$ is even. The first three come from the Ramanujan--Nagell equation $x^2+7=2^n$. In this situation $2$ splits in $K$, the element $\eta=(p^m+x\sqrt{-C})/2$ lies in $\OK$, and $(\eta)=\mathfrak p_2^{\,k}\mathfrak a^q$ with $\mathfrak p_2$ a prime above $2$ and $k\equiv-2\pmod q$; so the method of Lemma~\ref{replemma} breaks down (compare the remark after \cite[Lemma~3.1]{Patel}).
\end{remark}

\section{Examples}\label{sec:examples}

All six families are non-empty. In each of the following examples the hypotheses of Theorem~\ref{mainthm} hold; in particular $p\nmid C$, $\gcd(x,y)=1$ and $y$ is odd. The class numbers were computed by counting reduced binary quadratic forms of the discriminant of $K$.

\begin{center}
\renewcommand{\arraystretch}{1.25}
\begin{tabular}{@{}llllc@{}}
\toprule
Family & $C$ & $u$ & identity & $h(-C)$\\
\midrule
(I), $q=3$ & $5$ & $2$ & $5\cdot34^2+59^2=21^3$ & $2$\\
(I), $q=5$ & $5$ & $2$ & $5\cdot410^2+1801^2=21^5$ & $2$\\
(II) & $105$ & $4$ & $105\cdot53772^2+71^2=6721^3$ & $8$\\
(III) & $7$ & $2$ & $7\cdot430^2+3^6=109^3$ & $1$\\
(I$'$) & $11$ & $3$ & $11\cdot36^2+37^2=25^3$ & $1$\\
(II$'$) & $11$ & $1$ & $11\cdot8^2+5^2=9^3$ & $1$\\
(II$'$) & $19$ & $1$ & $19\cdot16^2+7^2=17^3$ & $1$\\
(III$'$), $\varepsilon=1$ & $11$ & $1$ & $11\cdot2^2+3^4=5^3$ & $1$\\
(III$'$), $\varepsilon=-1$ & $19$ & $1$ & $19\cdot28^2+3^6=25^3$ & $1$\\
(III$'$), $\varepsilon=1$ & $35$ & $1$ & $35\cdot26^2+3^6=29^3$ & $2$\\
\bottomrule
\end{tabular}
\end{center}

Proposition~\ref{congruence} can be checked on the generic rows: $\leg{-5}{59}=-1$ and $3\mid60$; $\leg{-5}{1801}=1$ and $5\mid1800$; $\leg{-11}{37}=1$ and $3\mid36$.

The converse part of Theorem~\ref{mainthm} also produces solutions when $q\mid h(-C)$, which are then outside the scope of the classification; for example $(C,u)=(61,2)$ and $(547,2)$ in family (III) give
\[
61\cdot3886^2+3^8=973^3,\qquad 547\cdot34990^2+3^{10}=8749^3,
\]
with $h(-61)=6$ and $h(-547)=3$.

\begin{remark}\label{computation}
As a computational check on Theorem~\ref{mainthm}, we enumerated all solutions of $Cx^2+p^{2m}=y^q$ with $\gcd(x,y)=1$, $y$ odd and $p\nmid C$, for odd square-free $5\le C\le119$ and $q=3,5,7,11,13$, with $y\le700,80,25,7,5$, respectively. For each $q$ and each odd $y$ in the indicated range, we enumerated the odd primes $p$ and integers $m\ge1$ with $p^{2m}<y^q$, tested in exact integer arithmetic whether $y^q-p^{2m}=Cx^2$ for an admissible $C$ and a positive integer $x$, and imposed $\gcd(x,y)=1$ and $p\nmid C$. The search found $208$ solutions. Of these, $58$ satisfy $q\nmid h(-C)$, and each of them belongs to exactly one of the six families ($35$ to (I), $1$ to (III), $14$ to (I$'$), $5$ to (II$'$) and $3$ to (III$'$); the smallest example of (II) lies outside the search range). The remaining $150$ have $q\mid h(-C)$, and some of these fit none of the families, for example $23\cdot4^2+19^2=9^3$ with $h(-23)=3$.
\end{remark}

\section{Concluding remarks}\label{sec:final}

\textbf{Composite exponents.} If $n$ is odd and $q\mid n$ is a prime with $q\nmid h(-C)$, then a solution of $Cx^2+p^{2m}=y^n$ gives the solution $(x,y^{n/q})$ of \eqref{eqthm}. Hence Theorem~\ref{mainthm} reduces the composite case to deciding when the values of $y$ in the six families, namely $Cu^2+1$, $4Cu^2+1$, $4Cu^2-3$, $(Cu^2+1)/4$, $Cu^2-2$ and $Cu^2\mp6$, are perfect powers subject to the accompanying conditions. This leads to further exponential equations, which we do not pursue here.

\textbf{Open questions.} It would be interesting to decide for which $C$ and $u$ the quantities
\[
|\I(C,u,1,q)|,\qquad 3Cu^2+1,\qquad (3Cu^2-1)/8,\qquad 3Cu^2-8
\]
are prime powers, to describe the solutions with $y$ even when $C\equiv7\pmod 8$, and to treat the case $q\mid h(-C)$.

\end{document}